\documentclass[11pt]{article}
 
 \usepackage{amsfonts,amssymb}

 \newcommand{\Z}{\mathbb{Z}}

 \newtheorem{theorem}{Theorem}[section]
 \newtheorem{corollary}[theorem]{Corollary}
 \newtheorem{lemma}[theorem]{Lemma}
 \newtheorem{proposition}[theorem]{Proposition}

 \newtheorem{definition}[theorem]{Definition}
 
 \newtheorem{example}[theorem]{Example}

 \newtheorem{remark}[theorem]{Remark}
 
 \newtheorem{conjecture}[theorem]{Conjecture}

 \input{epsf.sty}

\begin{document}

 \title{Virtual Knot Invariants from Group Biquandles and Their Cocycles
  }

 \author{
 J. Scott Carter\footnote{Supported in part by NSF Grants DMS \#0301095 and \#0603926.}
 \\University of South Alabama 
 \and 
 Mohamed Elhamdadi\footnote{Supported in part by 
 University of South Florida Faculty Development Grant} 
 \\ University of South Florida
 \and
 Masahico Saito\footnote{Supported in part by NSF Grants  DMS \#0301089 and \#0603876.}
 \\ University of South Florida
 \and
 Daniel S. Silver\footnote{Supported in part by NSF Grant DMS \#0304971}
 \\University of South Alabama
 \and
 Susan G. Williams\footnote{Supported in part by NSF Grant DMS  \#0304971}
 \\University of South Alabama
 }

 \maketitle

 \begin{abstract}

A group-theoretical method, via Wada's representations, is presented to distinguish
Kishino's virtual knot from the unknot. Biquandles are constructed for any group using
Wada's braid group representations.
Cocycle invariants for these biquandles are studied.
These invariants are applied to show the non-existence of Alexander numberings and checkerboard  colorings. 
 \end{abstract}

 \section{Introduction}

 The purposes of this paper include defining biquandle  structures on groups, and giving a group-theoretic proof that Kishino's virtual knot is non-trivial.
 A biquandle structure or a birack structure is related to 
 solutions to the 
 set-theoretic 
 Yang-Baxter equation (SYBE). 
 Given 
 an invertible solution that satisfies an additional condition (corresponding to a Reidemeister type I move), 
 we obtain a biquandle, and every biquandle gives a solution to the SYBE. Most examples that 
 were known up to this point came
 from generalizations of the Burau  representation. 
 The principal examples that we consider in this paper 
 come from Wada's representations of braid groups as free group automorphisms.  
 Our first example indicates  that one of these representations can be used to 
 distinguish Kishino's virtual knot from the unknot. 
 By abelianizing such groups, we recover an analog of Burau matrix. 
 Using these examples, we  construct and calculate cocycle invariants 
 that come from the homology theory of biquandles.
 As applications we give obstructions to checkerboard colorability and (mod 2)-Alexander numberings of arcs of virtual knots. 

 In Section~\ref{groupsec} we examine Wada's group invariants for virtual knots.
 The biquandle structures are defined on any group in Section~\ref{biqsec}
 using Wada's representations,
  and colorings of virtual knot diagrams by such biquandles are studied.
  Cocycle invariants are constructed and applied in Section~\ref{cocysec}.

 \section{Wada's group invariants for virtual knots} \label{groupsec}

 Wada~\cite{Wada}  
 considered representations of braid groups to the automorphism 
 groups of  free groups  of the following type. Fix finite words $u(x,y)$ and $v(x,y)$ in $x$ and $y$. Assume that the map $x \mapsto u(x,y)$, \ $y \mapsto v(x,y)$ is an automorphisms of the free group $F_2= \langle x,y \rangle$.
Wada's representation, $\rho:B_n\rightarrow {\rm Aut}(F_n)$  of the $n$-string braid group, $B_n$, to the automorphism 
group %%JSC
of the free group $F_n=\langle x_1, \ldots, x_n \rangle  $ is defined as follows.
 The standard  generator 
 $\sigma_i$ 
 is sent to the 
 automorphism
\begin{eqnarray*}
x_i \rho (\sigma_i) &=& u(x_i, x_{i+1}), \\
x_{i+1} \rho (\sigma_i) &=& v(x_i, x_{i+1}), \\
x_j \rho (\sigma_i) &=& x_j,  \quad ( j\neq i, i+1 ).
\end{eqnarray*}
 Such representations are also studied independently by A. J. Kelly~\cite{Kelly}
 as mentioned in the review article by Przytycki~\cite{Pr} of Wada's paper.

 A set of necessary conditions for the words $u(x,y)$ and $v(x,y)$ to define braid group representations was given by Wada. These are labeled $T$, $M$, $B$ below, to indicate the conditions that come from the top, middle, and bottom arcs, respectively,  in the type III Reidemeister  move:
 \begin{itemize}
 \item[T:] \hfil $u(u(x,y),u(v(x,y),z))= u(x, u(y,z))$ 
 \item[M:] \hfil $v(u(x,y),u(v(x,y),z)) = u(v(x,u(y,z)),v(y,z))$
 \item[B:] \hfil $v(v(x,y),z)=v(v(x,u(y,z)),v(y,z))$
 \end{itemize}
 It is interesting to compare these to the conditions for a birack
 \cite{CESbiq,KauRad,NV}. The conditions $T,M,B$ are called {\it Wada's conditions.}

\begin{figure}[htb]
\begin{center}
 \mbox{
\epsfxsize=4.5in
 \epsfbox{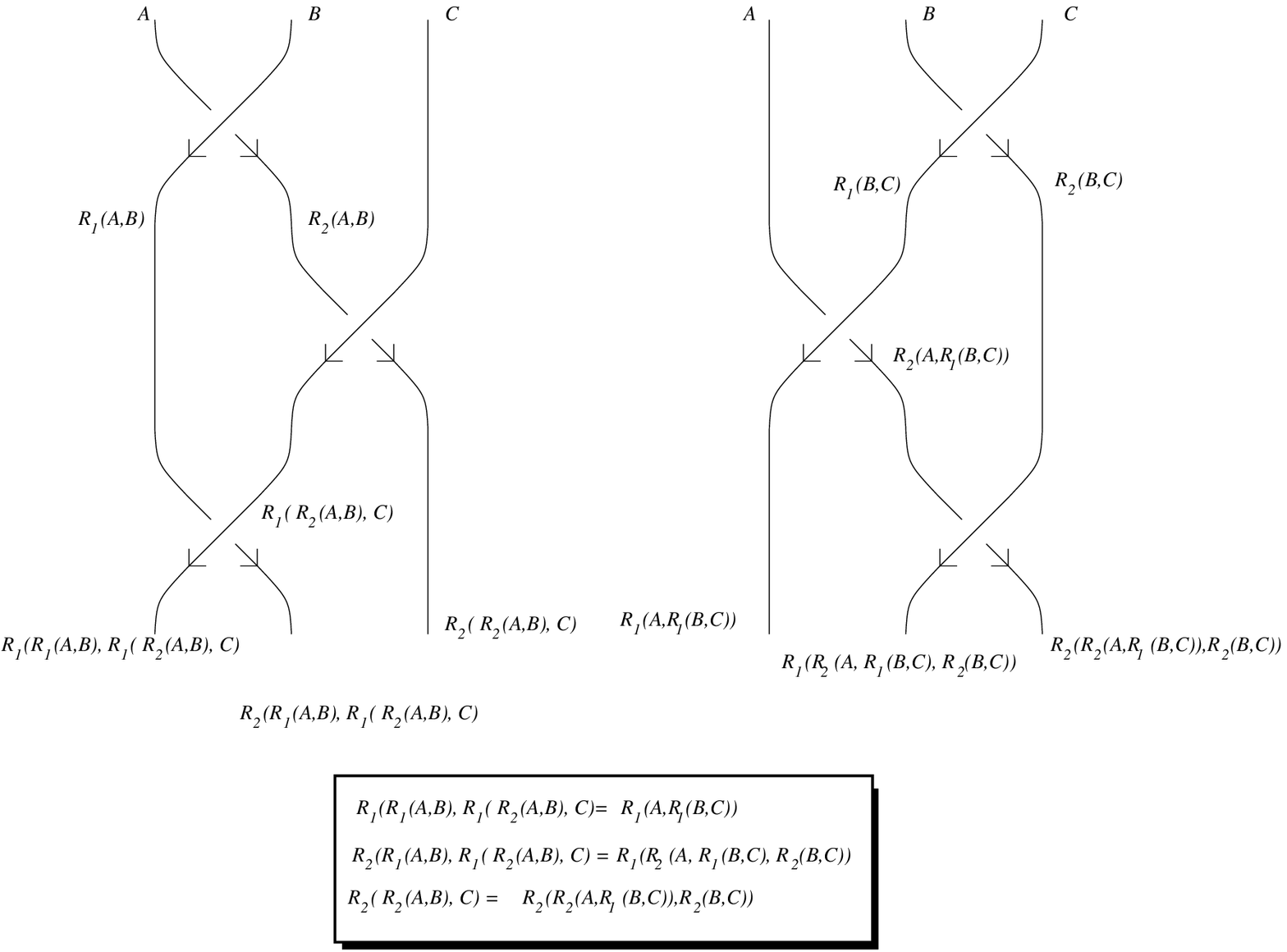} 
 }
\end{center}
\caption{The set-theoretic Yang-Baxter equation}
\label{wadarules} 
\end{figure}

Figure~\ref{wadarules} below indicates 
 the correspondence with the Reidemeister type III move
using the biquandle notation $(R_1(A,B),R_2(A,B))$ introduced in Section~\ref{biqsec}.

\begin{figure}[htb]
 \begin{center}
 \mbox{
 \epsfxsize=3in
 \epsfbox{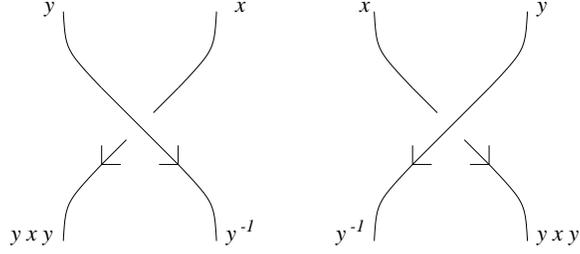} 
 }
 \end{center}
 \caption{Wada's braid representation, type $W_1$}
 \label{wadacross1} 
 \end{figure}

\begin{figure}[htb] 
 \begin{center}
 \mbox{
 \epsfxsize=3in
 \epsfbox{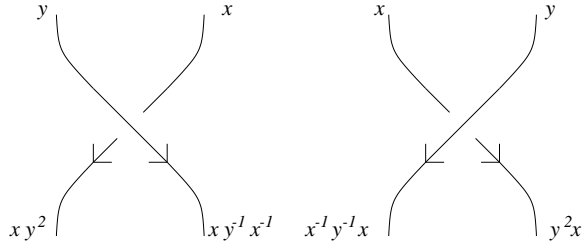} 
 }
 \end{center}
 \caption{Wada's braid representation, type $W_2$}
 \label{wadacross2} 
 \end{figure}

 The last two types of representation  
 that Wada considered (types 6  and 7 in his notation)  
 are depicted and named $W_1$ and $W_2$ in Figs.~\ref{wadacross1}
  and \ref{wadacross2}, respectively. 
  In the type $W_1$ representation $u(x,y)=y^{-1}$ 
 and $v(x,y)=yxy$ while in the $W_2$ representation $u(x,y)=x^{-1}y^{-1}x$ and $v(x,y)=y^2 x$. The left-hand side of each figure shows the words corresponding to the inverse automorphism of $F_2$. 
 Another braid group 
 representation that Wada considered 
 yields the {\it the core group}, 
 ${\mbox{\rm Core}}(K)$, of a knot $K$. 
 In this representation $u(x,y)=y$ and $v(x,y)=yx^{-1}y$. If $K$ is a classical knot 
 (one without virtual crossings, see the next paragraph), 
 then Kelly showed these three representations 
 give rise to the same knot invariants. 

 We generalize Kelly's result in Theorem~\ref{danandsusan} below.
 It is a pleasant exercise to show that the core group representation, $W_1,$ and $W_2$ satisfy Wada's conditions.

 Virtual knots and links were introduced by Kauffman (see \cite{KauB} for example). 
 A {\it virtual knot or link diagram} consists of a 4-valent 
 graph embedded in the plane with two types of vertices called {\it classical crossings} and {\it virtual crossings}. A classical crossing has two non-adjacent edges at a vertex being distinguished as an {\it over-crossing arc} with the remaining two forming a pair of {\it under-crossing arcs}. The under-crossing arcs are indicated in a drawing by not connecting them to the vertex. In a drawing of a  virtual crossing, the vertex is encircled. Two non-adjacent edges are 
 regarded as 
 a single edge. Thus the {\it edges} in a virtual diagram start and end at classical crossings and 
 may pass through a virtual crossing.
 A {\it virtual knot} or {\it virtual link} is an equivalence class of  virtual diagrams modulo the equivalence relation generated by the virtual Reidemeister moves.

 Let ${\cal E}$  
 denote the set of edges of a given diagram.
 Then Wada's group invariants are defined by assigning 
 generators to elements of ${\mathcal E}$ and assigning 
 relations to crossings. 
 Specifically, let ${\mathcal E}=\{ x_1, \ldots, x_m\}$.
 These letters will be  also used as generators.
At 
a positive crossing as depicted in the 
right of Fig.~\ref{wadacross1} or Fig.~\ref{wadacross2}
 denote the (generators corresponding to)    top 
 edges 
 $x, y \in {\mathcal E}$ as depicted. Let $u,v \in {\mathcal E}$ be the 
 bottom 
 edges. 
 Then the relations $r_j, r_j'$ corresponding to this crossing
 are $
 u=u(x,y)$, $
 v=v(x,y)$, 
respectively, where 
$u(x,y)$ and 
$v(x,y)$ 
 are words 
 in $x,y$ that define Wada's representation. 
 For example, for type $W_1$, 
 the relations are 
  $u=y^{-1}$ and $v=yxy$, cf.  
the right of 
Fig.~\ref{wadacross1}. 
At negative crossings,  the words on the left are used. 
 The following is  a straightforward calculation.

 \begin{lemma} 
 Wada's group invariants are well-defined for virtual knots 
and links. 
 \end{lemma}

 Denote by $W_i(K)$, $i=1,2$, the group invariants that are defined by using 
 the type $W_i$ represention. 
 Note that these groups are free groups of rank $n$ for the unlink
 of $n$ components for any positive integer $n$.

 \begin{figure}[htb]
 \begin{center}
 \mbox{
 \epsfxsize=3.5in
 \epsfbox{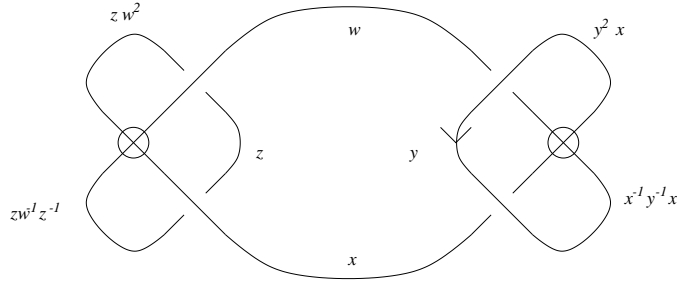} 
 }
 \end{center}
 \caption{A group invariant for Kishino's knot}
 \label{kishinogpnu} 
 \end{figure}

 As an example we compute $W_2(K)$ for 
 a virtual knot $K$ widely known as Kishino's knot, depicted in 
 Fig.~\ref{kishinogpnu}.
 It has  attracted attention for its 
 remarkable property that it is a connected sum of two diagrams 
 of the trivial knot; it has  trivial Jones polynomial, 
 trivial fundamental group and quandle.
 Furthermore, it  can be seen by computing presentations  
 that both 
 its core group and its Wada group of type $W_1$ are infinite cyclic.

 Kishino's knot 
 was 
 first proved to be non-trivial  in \cite{KS} by using 
 the Jones polynomial of the 3-strand parallel of the virtual knot. 
 Later other methods were discovered, including a computation of the minimal genus of the surface 
 on which 
 the virtual knot lies~\cite{Kado}, biquandles~\cite{BF}, 
 and 
 the bracket on surfaces~\cite{Heather}. 
 Creel and Nelson \cite{CrN} used Yang-Baxter cocycles \cite{CESbiq} 
 to distinguish Kishino's example from the unknot.
 See also \cite{DyeKauf}.

 \begin{proposition}\label{kishinoprop}
 For the Kishino knot $K$, the group invariant $W_2(K)$ is not 
 infinite cyclic.
 \end{proposition}
 {\it Proof.\/}
 In Fig.~\ref{kishinogpnu}, Kishino's knot is depicted with labels on 
 four of the edges  induced by the $W_2$ representation. The first two relations below are induced at the top right crossing, and the last two are induced at the bottom left crossing:
 \begin{eqnarray}
 w&=&(x^{-1}y^{-1}x)(y^2x)^2, \\
 y&=&(x^{-1}y^{-1}x)(y^2x)^{-1}(x^{-1}y^{-1}x)^{-1},\\
 x&=& (zw^{-1}z^{-1})^{-1}(zw^2)^{-1}(zw^{-1}z^{-1}),\\
 z&=&(zw^2)^2(zw^{-1}z^{-1}).
 \end{eqnarray}
 These relations are induced from the following considerations.
 The edge with label $w$  (on the upper right of the diagram) is the target of the under-crossing edge, and  
 the arc labeled $y$ is the target  of the over-crossing edge. On the upper right crossing the edges are labeled $y^2 x$ and $x^{-1}y^{-1}x$ since these are the target edges of the lower right crossing.
 These give the first two relations above. 
 Similar considerations hold on the lower left for the arcs with labels $z$ and $x$.

 The relation (4) simplifies to  $(4'):$ $zw=(w^2z)^2$.
 Using this, (3) simplifies to  $(3'):$ $x=(w^2 z)^{-1}$, or $z=(xw^2)^{-1}$, which 
 simplifies $(4')$ to  $(4''):$ $xw^2=wx^2$.
 The relation (2) simplifies to  $(2'):$ $x^{-1}yx=y^3xy$.
 Repeated use of $(2')$ simplifies (1) to  $(1'):$ $w=y^{-3} x^{-1} y^{-3}$. 
 Hence we obtain 
 $$W_2(K)=\langle \  
  x,y\ | \   x^{-1}yx=y^3xy, \  x( y^{-3} x^{-1} y^{-3}) ^2=(y^{-3} x^{-1} y^{-3})x^2\ \rangle.$$
 Under the substitution $a=xy^3$, the first relation becomes
  $ya=a^2y$, and the second relation is rewritten as 
 $ y^3 a^{-1} y^3 a^2= a^2 y^3 a y^3 $ (take the inverse of both sides of 
 the relation and simplify after substitution).
 The first relation can be seen as a commutation relation,
 and in general $ya^m=a^{2m}y$ for any positive integer $m$.
 It follows that $y^3 a^{15}=a^{10} y^3$.
 Hence the group has presentation 
 $$W_2(K)=\langle \  
 a, y \ | \ ya=a^2y, \ y^3 a^{15}=a^{10} y^3\ \rangle.$$

The first relation implies that $y^3 a y^{-3} = a^8$. Hence the second relation can be replaced by 
$a^{120}y^3=a^{10}y^3$
which is equivalent to 
 $a^{110} =1$. Using the first relation again, we can write $a^{110} = (a^{55})^2 = y a^{55} y^{-1}$. Therefore we can replace $a^{110}=1$ by the conjugate relation $a^{55}=1$.
Thus
$$W_2(K) = \langle a, y \mid yay^{-1} = a^2, a^{55}\rangle.$$ 
Regarded additively, this group is simply the semidirect product of ${\mathbb Z}_{55}$  by an infinite cyclic group $\langle y \mid \rangle$, the action of $y$ being multiplication by $2$. In particular, 
$W_2(K)$ is not cyclic. $\Box$

 \bigskip 

 It is known~\cite{Wada,Kelly}  that if $K$ is classical, then 
 the groups $W_i(K)$, $i=1,2$, are isomorphic to the core group 
 ${\rm Core}(K)$. 
 Kelly's result can be extended in some cases to the virtual category. In this way a necessary criterion for Alexander numbering is formulated.

\begin{figure}[htb]
 \begin{center}
 \mbox{
 \epsfxsize=2.5in
 \epsfbox{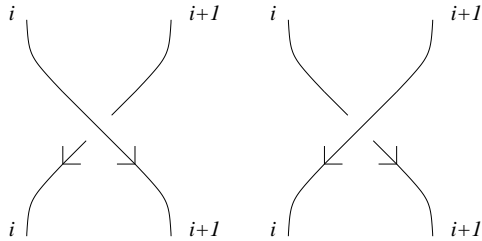} 
 }
 \end{center}
 \caption{Alexander numbering}
 \label{LXN} 
 \end{figure}

 \begin{definition}
 {\rm  \cite{DanSusan} An {\it Alexander numbering} of 
 edges 
   ${\mathcal E}$ of an
 oriented  virtual link diagram
 is an assignment of integers on 
 edges 
 such that at every crossing the following condition 
 is satisfied. When a (positive or negative) crossing is placed so that 
 both 
 edges are oriented downward, the left 
 edges 
 receive the same integer $i$,
 and then the right 
 edges receive the integer $i+1$ (see Fig.~\ref{LXN}).

 The numbering of an edge remains unchanged when it passes through a virtual crossing. 
 A {\it (mod 2)-Alexander numbering} is defined similarly, but the indices on the edges are integers modulo $2$. 
 } \end{definition}

 \begin{definition}{\rm  \cite{DanSusan} 
  An oriented virtual link is {\it almost
 classical} if it has a diagram that admits an Alexander numbering.
It is {\it almost classical (mod 2)} if it admits a (mod 2)-Alexander numbering. 
 } \end{definition}

 An example is given in 
 \cite{DanSusan} 
 of a 
 virtual knot 
 that is almost classical but is shown 
%%%by 
 by Alexander groups 
 to be non-classical.

 \begin{proposition} \label{danandsusan}
 If an oriented diagram 
for  a virtual knot 
 $K$ has a (mod 2)-Alexander numbering, then the group 
 $W_1(K)$, 
 is isomorphic to the core group of $K$. \end{proposition}
%%%close %%%and \it
{\it Proof.} We can write a presentation for $W_1(K)$ with a generator for each edge of the diagram and a pair of relations at each crossing. These are of the form $c=u(a,b)$, $d=v(a,b)$, with $u(a,b)=b^{-1}$ and $v(a,b)=bab$. If the crossing is positive then $a$ and $b$ are the generators corresponding to the left and right top edges, and $c$ and $d$ correspond to the bottom left and right edges, respectively. For a negative crossing, $a$ and $b$ correspond to the bottom left and right edges, and $c$ and $d$ are the left and right top edges, respectively. The core group has a similar presentation, but there $u$ and $v$ are replaced by $u'(a,b)=b$ and $v'(a,b)=ba^{-1}b$. 

 Given a (mod $2$)-Alexander numbering for the diagram $K$, we obtain an isomorphism from $W_1(K)$ to the core group by sending each generator $x$ to $x^{-{\rm ind}}$ where ind denotes the index on the corresponding edge. If $a$ and $c$ (which are on the left) have index $0$, then the crossing relations are sent to $c=u(a,b^{-1})=u'(a,b)$ and $d^{-1}=v(a,b^{-1})=(v'(a,b))^{-1}$. Similarly if $a$ and $b$ have index $1$, then the relations are sent to $c^{-1}=u(a^{-1},b)=u'(a,b))^{-1}$ and $d=v(a^{-1},b)=v'(a,b)$. Thus in either case the image of the pair of $W_1$ relations at a crossing is the pair of core relations at the same crossing, up to a trivial rewriting. The (mod 2)-Alexander numbering allows us to make a global choice for the images of the generators. $\Box$

 \begin{example} {\rm If a virtual knot or link $K$ does not have a diagram admitting a (mod $2$)-Alexander numbering, then $W_1(K)$ need not be isomorphic to the core group of $K$.

 \begin{figure}[htb]
 \begin{center}
 \mbox{
 \epsfxsize=2.5in
 \epsfbox{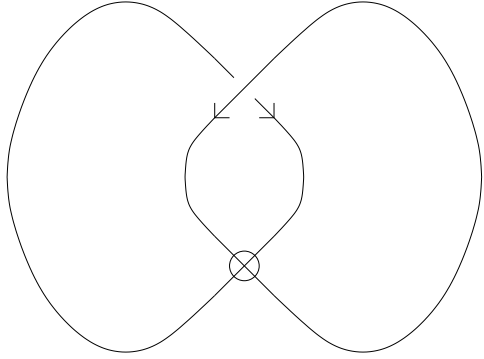} 
 }
 \end{center}
 \caption{A virtual link for which $W_1(K)$ and ${\rm Core}(K)$ are distinct}
 \label{vhopflink} 
 \end{figure}

 The simplest example is provided by the virtual Hopf link in Fig.~\ref{vhopflink}. 
An easy computation shows that $W_1(K)$ is isomorphic to $\langle x, y \mid xy = yx, y^2\rangle \cong {\mathbb Z} \oplus {\mathbb Z}_2.$ On the other hand,  the core group of $K$ is isomorphic to $\langle x, y \mid x y^{-1} = y x^{-1} \rangle.$  Letting $z = yx^{-1}$, we immediately see that the latter group is  isomorphic to $\langle x, z\mid  z^2 \rangle$,  the free product of ${\mathbb Z}$ and ${\mathbb Z}_2$.  

An example for which $K$ is a 
virtual 
knot is given here.
Consider the knot $K$ in Fig.~\ref{vgranny}, a connected sum of two virtual trefoils. A straightforward calculation, similar to that above, shows that $W_1(K)$ 
is the classical trefoil knot group $\langle x, y \mid xyx=yxy\rangle$ while the core group ${\rm Core}(K)$ is the free product of ${\mathbb Z}$ and ${\mathbb Z}_3$.

 \begin{figure}[htb]
 \begin{center}
 \mbox{
 \epsfxsize=2.5in
 \epsfbox{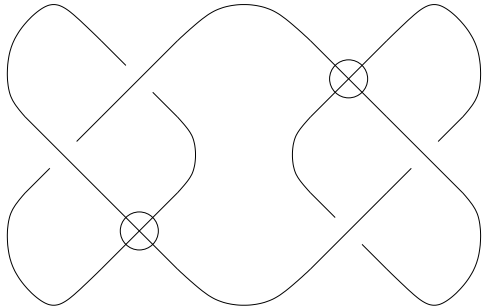} 
 }
 \end{center}
 \caption{A virtual knot for which $W_1(K)$ and ${\rm Core}(K)$ are distinct}
 \label{vgranny} 
 \end{figure}

}\end{example}

 \bigskip

 We conclude this section with:

 \begin{conjecture}\label{Zconj}
 The free part of the abelianization of 
a
Wada group 
is $\Z $.
 \end{conjecture}

 \bigskip

 \section{Biquandles from Wada's representations}\label{biqsec}

 In this section we define biquandle structures on any group using
 Wada's representations and study colorings  of virtual knots by such biquandles.

 A {\it quandle}, $X$, is a set with a binary operation $(a, b) \mapsto a * b$
 such that

 (I) For any $a \in X$,
 $a* a =a$.

 (II) For any $a,b \in X$, there is a unique $c \in X$ such that 
 $a= c*b$.

 (III) 
 For any $a,b,c \in X$, we have
 $ (a*b)*c=(a*c)*(b*c). $

 A {\it rack} is a set with a binary operation that satisfies 
 (II) and (III).
 Racks and quandles have been studied in, for example, 
 \cite{Brieskorn,FR,Joyce,K&P,Matveev}.
Generalizations of racks and quandles 
have been studied in several papers. 
 Here we follow descriptions in \cite{BF,FJKW}.

 \begin{definition} {\bf  \cite{BF,FJKW}} {\rm 
 A 
{\it birack} is a 
set $X$ together with a mapping 
 $R: X \times X \rightarrow X \times X$ 
that has the following properties.

 \begin{enumerate}
 \setlength{\itemsep}{-5pt}
 \item
 The map $R$ is invertible.
 The inverse of $R$ is denoted by $\bar{R}: X \times X \rightarrow X\times X$. 
We will write 
\\[-3mm]
 $$R(A_1, A_2)=(R_1 (A_1, A_2), R_2 (A_1, A_2) )= (A_3, A_4) , $$
 where $A_i \in X$ for $i=1,2,3,4$.

 \item 
 For any $A_1, A_3 \in X$ there is a unique $A_2 \in X$ such that
 $R_1 (A_1, A_2) = A_3$.
 We say that $R_1$ is {\it left-invertible. }

 \item 
 For any $A_2, A_4 \in X$ there is a unique
  $A_1 \in X$ such that $R_2 (A_1, A_2)= A_4$. 
 We say that  $R_2$ is {\it  right-invertible.} %%JSC space

 \item 
 $R$ satisfies the set-theoretic Yang-Baxter equation:\\[-3mm]
 $$ (R \times 1) (1 \times R) (R \times 1)
 =  (1 \times R)(R \times 1) (1 \times R) ,$$
 where $1$ denotes the identity mapping.
 \end{enumerate}

 } \end{definition}

 \begin{definition}{\bf \cite{BF,FJKW} } \label{biqdef}
 {\rm 
 A biquandle $(X, R)$ is a birack with the following property,
 called the {\em type I condition}:
 For any $a \in X$ there are unique elements $x_a$ and $y_a$ such that $x_a=R_1(x_a,a)$, $a=R_2(x_a,a),$  $y_a=R_2(a,y_a)$, 
and $a=R_1(a, y_a)$. 
 } \end{definition}

 \begin{figure}[htb]
 \begin{center}
 \mbox{
 \epsfxsize=3.5in
 \epsfbox{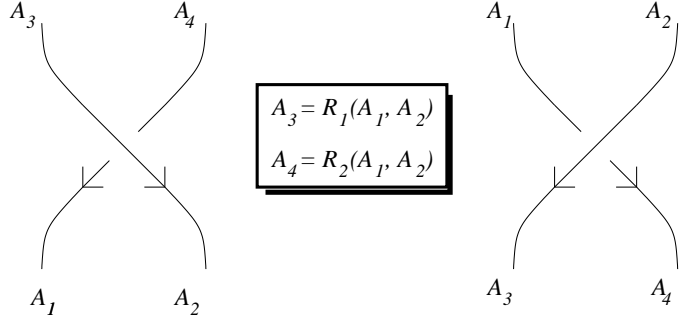} 
 }
 \end{center}
 \caption{A coloring by birack elements }
 \label{Rmatrix} 
 \end{figure}

Direct calculations imply the following. 

 \begin{lemma} 
Suppose that $u(x,y)$ and $v(x,y)$ are words in $F_2$ such that $x \mapsto u(x,y)$ and $y\mapsto v(x,y)$ define isomorphisms of the free group 
that 
satisfy Wada's conditions. Then for any group $G$, the binary operations
 $R_i: G \times G \rightarrow G$, $i=1,2$, 
given  by 
$$R_1(x,y)=u(x,y), \quad R_2(x,y)=v(x,y)$$
define 
a birack structure on $G$. 
The biracks corresponding to the type $W_1$ and $W_2$ representations,
 $$ R_1(x,y)=y^{-1}, \quad R_2(x,y)=yxy, $$  
 $$\mbox{\rm or} \ \ R_1(x,y)=x^{-1} y^{-1} x, \quad R_2(x,y)=y^2 x, $$  
respectively, are biquandles. 
 \end{lemma}  

 Call the first and second, a group  biquandle of type  $W_1$ and $W_2$, respectively.
 Note that if the group $G$ is abelian, then both types define the same 
 map $R(x,y)=(-y, 2y+x)$. We call this an {\it abelian Wada biquandle.}

  Let  $(X, R)$ be a biquandle and ${\mathcal E}$ be the set of edges of
  a virtual knot diagram $K$ (an edge runs from classical crossing to classical crossing).

 \begin{definition}{\bf \cite{BF,FJKW} } \label{biqcolor} 
 {\rm 
 A coloring ${\mathcal C}:  {\mathcal E} \rightarrow X$ is a map such
  that at every crossing, the conditions depicted in Fig.~\ref{Rmatrix}
   are satisfied where $(A_3, A_4)= R(A_1,A_2),$  and $A_1,\ldots, A_4 \in X.$
  } \end{definition}

 Let ${\rm Col}_X(K)$ denote the set of colorings of a virtual knot diagram
 $K$ by a biquandle $X$. 
If two diagrams are related by a generalized Reidemeister move, there is a one-to-one correspondence between their sets of colorings, so that  
 the 
cardinality
of colorings $ | {\rm Col}_X(K) |$ 
  is an invariant of virtual knots 
 (see \cite{BF,FJKW}, for example).

\begin{remark} 
{\rm If the Wada group invariant 
of $K$  
associated to $R$ 
is $\Gamma$, then ${\rm Col}_X(K)$ can be identified with ${\rm Hom}(\Gamma, G)$.
%%the set of homomorphisms of biquandles appropriately defined. %% 
%%NO, it is true at the group level. JSC

}\end{remark}

 \begin{definition} {\rm
 A  
coloring of a knot $K$ by the biquandle (G,R) with $G=\Z_n$  or $G=\Z$ --- in this case $R(x,y)=(-y,2y+x)$ ---  is called a {\it cyclic coloring.}} \end{definition}

 Note that any virtual diagram has a cyclic coloring by $\Z_n$, for any $n$,
 such that  every arc receives the color $0$. 
 On the other hand, if a diagram has a cyclic coloring with a single color 
 $m\neq 0 \in \Z_n$,
 then either the diagram is the  diagram  
 of the unknot without a real crossing, or 
 $n=2m$, 
 as
 $m=x=y=-y=x+2y$ must be satisfied at every real crossing.

 \begin{proposition}
 Any virtual link has a coloring by the infinite cyclic abelian Wada biquandle $\Z$
 such that at least one color is not $0$.
 \end{proposition}
 {\it Proof.\/} 
 Any virtual link is represented as the closure of a virtual braid~\cite{Seiichi}
 of $n$ strings for some positive $n$,
 that is, a sequence of ordinary braid generators  and 
 virtual crossings. Denote by $\sigma_i$ and $v_i$ the standard braid generator 
 and the virtual crossing, respectively, at the $i\/$th and $(i+1)\/$st strings,
 $i=1, \ldots, n-1$. Let $\beta$ be an $n$-string virtual braid  
 whose closure represents a given virtual link.
 Colors  at the top strings are  represented by a vector $w \in \Z^n$,
 and subsequent colors are computed by matrices.
 The matrix  corresponding to $\sigma_i$ is 
$\left[ 0  \; -1 \atop  1 \; \; \  2 \right]$ 
placed at the $(i, i+1)$ block
 of the identity matrix of size $n$, and for $v_i$, the transposition matrix
$\left[  0 \ 1\atop 1 \ 0  \right]$
is placed.
 Hence the vector $w'$ representing the colors of the bottom strings
 of the braid is computed by $Mw$, where $M$ is the matrix 
 representing $\beta$, which is a product of the above described matrices.
 Note that each matrix above has the property that $w_1=[1, \ldots, 1]$
 (the row vector with every entry $1$) is a left eigenvector of eigenvalue $1$,
 so that $M$ has the same property. Since $M$ has the eigenvalue 
 $1$, there
 is a right eigenvector 
 of $M$ 
 over the field of rational numbers,
 with eigenvalue $1$.
 Multiplying by a common multiple of the denominators of the entries
 of such an eigenvector, we obtain an integral eigenvector. Using the entries to color the top strands of $\beta$, we obtain a coloring that extends to the closure. 
 $\Box$

 \begin{corollary}
 For any oriented virtual link $L$ and for $i=1,2$, 
 there are epimorphisms from the Wada groups 
 $G_i(L) \rightarrow \Z$ for $i=1,2$.
 \end{corollary}

 In \cite{Naoko} checkerboard colorable virtual knots were defined,
 and the properties of their Jones polynomials 
were studied. 
In \cite{KK} an abstract link diagram is constructed from a virtual diagram as follows. A surface with boundary, $F$,  is constructed as a handlebody. 
The $0$-handles correspond to the classical crossings in the virtual diagrams; the $1$-handles correspond to the edges in the diagram and four such $1$-handles are attached in the natural fashion to a $0$-handle as the  edges approach the vertices. The resulting diagram represents a link in the $3$-manifold that is  $F \times I$. 
 A virtual link diagram is {\it checkerboard colorable} if  the graph of the virtual link in $F$ is checkerboard colorable. Equivalently, the $1$-cycle represented by the link is null-homologous modulo 2.
If a virtual link possesses such a diagram, it is called
 {\it checkerboard colorable}. 
 The relation between checkerboard colorability and Alexander numberings 
  is 
 as follows.
 A virtual link  is checkerboard colorable if and only it has 
 an oriented   diagram  
that  admits 
 a (mod $2$)-Alexander numbering. For we can color 
the regions to the left of  the edges
that are 
labeled 
 $0$ by black, right by white, and  
use the opposite colors for those edges that are labeled $1$.

 \begin{lemma}\label{Foxlemma}
 If  virtual link $L$ is checkerboard colorable, then there is a one-to-one 
 correspondence between the set of colorings of any given diagram 
 of $L$ 
 by an  abelian Wada biquandle $\Z_n$ and the set of Fox $n$-colorings.
 \end{lemma}
 {\it Proof.\/}
 Let $L$ be an oriented virtual knot diagram with a 
(mod $2$)-Alexander numbering,
 so that every edge $\alpha \in {\mathcal E}$ is assigned $0$ or $1$, 
 denoted by $\epsilon(\alpha)$.
 Let ${\mathcal C}$ be a cyclic  coloring
 by an  abelian Wada biquandle $\Z_n$. 
 Define a Fox coloring ${\mathcal C}'(\alpha)$, $\alpha \in {\mathcal A}$, 
 by  ${\mathcal C}'(\alpha)=(-1)^{\epsilon(\alpha)}{\mathcal C}(\alpha)$.
 Then it is checked that this provides a one-to-one 
 correspondence.
 $\Box$

 We remark that this change of basis is the abelian version of 
 Proposition~\ref{danandsusan}.

 \begin{corollary}\label{Foxcor}
 If a virtual knot $K$ is checkerboard colorable, then 
 $K$ is Fox $n$-colorable (non-trivially) if and only if $| {\rm Col}_X(K) |>  n$,
 where $X=\Z_n$ denotes a cyclic  abelian Wada biquandle.
 \end{corollary}

 \begin{remark}\label{spanrem}{\rm
 If an oriented  diagram of a virtual link is checkerboard colorable, 
 then there is a coloring by the cyclic 
Wada biquandle $\Z$ 
 using only $\pm 1$. Let $L$ be a (mod 2)-Alexander numbering, 
 and simply assign $1$ for an arc $\alpha$ with $L(\alpha)=0$ and $-1$ 
 if  $L(\alpha)=1$. Thus there is a non-zero coloring by $\Z$ such that 
 the span (the largest integer minus the smallest integer that appear in the colors)
 of the coloring is $2$.

 For example, from Fig.~\ref{kishinogpnu}, by computing the abelianized elements, we 
 see that this diagram of Kishino's knot  
does not have 
 a non-zero coloring by $\Z$ whose span is $2$.
 The smallest span is $4$. This does not, however, prove that Kishino's knot is
 not checkerboard colorable, as there might be another diagram with span $2$.

 The {\it span} of an oriented virtual link is
 the minimal span among all diagrams and all non-zero colorings by $\Z$. We conjecture that the span of Kishino's knot is $4$.
 } \end{remark}

 \section{Cocycle invariants}\label{cocysec}

 A homology theory and cocycle invariants
 for biquandles were defined in \cite{CESbiq}.
 In this section, we study constructions of cocycles  for group biquandles.
 For 
the 
 purposes of this paper, we review definitions of cocycles only  as 
 functional equations, instead of going through homology theories,
in order to use them to define cocycle invariants. 

 Let $(X,R)$ be a biquandle and $A$ be an abelian group.
 A function $g: X  \rightarrow A$
 is called 
 a {\it Yang-Baxter} $1$-cocycle if it satisfies:
 \begin{eqnarray*}
 ( \delta_1 g) (x,y) & := &  g(x) + g(y) - g(R_1(x,y)) - g (R_2(x,y)) \; = \; 0. \\
 \end{eqnarray*}

 A function $f: X^2  \rightarrow A$ 
 is called 
 a {\it Yang-Baxter} $2$-cocycle if it satisfies:
 \begin{eqnarray*}
 (\delta_2 f ) (x,y,z)
 &:=& f(x,y)+f(R_2(x,y),z)+f(R_1(x,y), R_1( R_2(x,y),z) ) \\
  &- & \{ f(y,z) + f(x, R_1(y,z)) + f(R_2(x, R_1(y,z)), R_2(y,z))\}  \; = \; 0.
 \end{eqnarray*}

The $2$-cocycle condition corresponds to the Reidemeister type III move. See Fig.~\ref{wadarules}. The ordered arguments of a function, $f$, are the left and right in-coming labels on the arcs near a crossing $A,B,R_2(A,B), C$ {\it etc.}. The three positive terms in the cocycle condition, then, correspond to the crossings in the 
left-hand side 
of the type III move, and the negative terms correspond to the crossings on the 
right-hand side. 
The type II move is handled by assigning the inverse value of $f$ at negative crossings.
To define 
knot invariants, the following condition becomes necessary
 for invariance under the type I Reidemeister move.

 \begin{definition}{\bf \cite{CESbiq}}
  {\rm 
 Let $(X, R)$ be a biquandle, 
 and $A$ be an abelian group. 
 Recall 
 that for any $a \in X$ there are unique
 elements $x_a$ and $y_a$ such that $x_a=R_1(x_a,a)$, $a=R_2(x_a,a),$  $y_a=R_2(a,y_a)$, and $a=R_1(a, y_a)$.

 A Yang-Baxter $2$-cocycle $f $ is  
  said to satisfy the {\it type I condition} if 
 $f(x_a,a)=0$ and $f(a,y_a)=0$
  for every $a \in X$.
 } \end{definition}

 Cocycles of Alexander quandles that are in polynomial form
were  
considered by Mochizuki~\cite{Mochi} and have been used 
to define 
cocycle knot invariants. These have been shown to have a wide range of applications.
 Let $G=\Z_n$ or $\Z$. 
(The case $G=\Z$ is referred to as the case $n=0$.)
Recall that the abelian Wada biquandle is
 given by $R(x,y)=(-y,2y+x)$. 
 Then the $2$-cocycle condition 
  with the coefficient group $A=\Z_n$ reads 
 $$f(x,y)+f(2y+x,z)+f( -y, -z)=f(y,z) + f(x, -z)+f(x-2z,2z+y).$$
 In this case, the type I condition for $2$-cocycles becomes
 $f(x,-x)  =  0$.

 \begin{lemma} 
 For any non-negative integer $n$, the function
 $f(x,y)=x+y$ is a $2$-cocycle of  the 
abelian biquandle 
$G=\Z_n$ with values in the coefficient 
group 
$A=\Z_n$. 
 It  
satisfies the type I condition and  is not a coboundary that is, there is no function $g:G\rightarrow A$ such that $\delta_1 g=f$. 
 \end{lemma}
 {\it Proof.} 
 It is   a direct calculation to check the $2$-cocycle condition and the type I 
 condition formulated above. 
To determine that $f$ is not a coboundary, we find a $2$-chain that is a cycle and upon which $f$ evaluates non-trivially. The set of $2$-chains is the free abelian group generated by pairs of elements in the biquandle $X$. The boundary of a generating chain is given as $\partial_2 (x,y)= \{x\}+\{y\} - \{R_1(x,y)\} - \{R_2(x,y)\}$ --- the braces indicate that  the elements of $X$ are considered as generators of a free abelian group on $X$. In the current context, $\partial_2(x,y)=  \{x\}+\{y\} - \{-y\} - \{2y+x\}.$

 Consider the $2$-chain $c=(1,0) \in \Z_n \times \Z_n$. 
 Then $c$ is a $2$-cycle,, and $f(c)=1\neq 0$, so that $f$ is not a coboundary.
 $\Box$

 \begin{corollary}
 The two-dimensional homology group $H_2^{\rm YB}(G;A)$ 
of  the abelian Wada 
 biquandle $G=\Z_n$ with the coefficient group $A=\Z_n$ is non-trivial for any non-negative integer $n$.
 \end{corollary}

 \begin{remark} {\rm
Let $p$ be an odd prime and consider 
cocycles of 
the abelian biquandle $G=\Z_p$ 
that take values in $A=\Z_p$. 
Consider the expression $h(x,y)=(1/p)[(x^p+2y^p)-(x+2y)^p]$ 
that is inspired by Mochizuki's cocycle. Since the numerator 
is divisible by $p$, $h(x,y)$ takes integral values, which are then reduced modulo $p$. 
A direct calculation gives that $h$ is a $2$-cocycle. 
For $p=3$ and $5$, respectively, $h$ evaluates non-trivially
 on the $2$-cycles $(1,1)+2(2,2)$, $(1,1)+2(2,2)+4(3,3)$. 
From the proof of the above lemma, we see that the cocycles
$f(x,y)=x+y$ and $h(x,y)=(1/p)[(x^p + 2y^p)-(x+2y)^p]$
are linearly independent and not coboundaries.
Hence rank$(H_2^{\rm YB}(G; A))$ is at least two for $p=3, 5$.
 } \end{remark}

%%%%delete and replace%%%%%
 %The biquandle $2$-cocycle invariants were defined in \cite{CESbiq}
% as follows. 
% Let $K$ be a classical or virtual  knot or link diagram. Let a finite 
% biquandle $(X,R)$, and a $2$-cocycle 
% $\phi : X^2 \rightarrow A$, that satisfies the type I condition,
% be given,
% where $A$ is an abelian group.
%For for the rest of  
 %the paper, we take a multiplicative cyclic abelian group for $A$ 
%generated by $u$;  the elements of $A$ are powers of $u$. 
% Let  ${\cal  C}$ 
% denote a coloring ${\cal  C}: {\mathcal E} \rightarrow X$,
% where $ {\mathcal E} $ denotes the set of arcs of $K$. 
%%%%%%%%%%%

%%%%replace%%%%%
We recall the definition of {\it biquandle $2$-cocycle invariants} from \cite{CESbiq}.
 Let $K$ be a classical or virtual  knot or link diagram. Let 
 $(X,R)$ denote a finite biquandle, and let   $\phi : X^2 \rightarrow A$ denote $2$-cocycle 
 that satisfies the type I condition where $A$ is an abelian group which is written multiplicatively.
 Let  ${\cal  C}$ 
 denote a coloring ${\cal  C}: {\mathcal E} \rightarrow X$,
 where $ {\mathcal E} $ denotes the set of arcs of $K$. 
%%%%%%%%%%%
%%%close-up 
 For a positive crossing $\tau$ as depicted in the right of Fig.~\ref{Rmatrix},  
 let $\alpha$ and  $\beta$ be 
 the top left and right arcs labeled by $A_1$ and $A_2$, respectively, and 
 ${\cal C}(\alpha)=A_1$ and  ${\cal C}(\beta)=A_2$ as depicted.
 For a negative crossing at the left of  Fig.~\ref{Rmatrix},
 let $\alpha$ and  $\beta$ be 
 the bottom left and right arcs, also  labeled by $A_1$ and $A_2$ respectively,
 so that  the coloring is also given by 
 ${\cal C}(\alpha)=A_1$ and  ${\cal C}(\beta)=A_2$.
%%Denote  the value of the cocycle $\phi,$ by $\phi(A_1, A_2)=u^{m} %%\in A$ for some positive integer $m$. 
 A {\it (Boltzmann) weight} $B(\tau, {\cal C})$ 
 (that depends on $f$)
 at a  crossing $\tau$ is defined 
 by $B(\tau, {\cal C}) =u^{{\epsilon (\tau) m}}$,   where $\epsilon (\tau)=1$ or $-1$ if $\tau$ is positive or negative, 
 respectively.

 The {\it (Yang-Baxter) cocycle knot invariant} is defined by 
 the state sum expression  
 $$
 \Phi_{\rm YB} (K) = \sum_{{\cal C}}  \prod_{\tau}  B( \tau, {\cal C}). 
 $$
 The product is taken over all crossings of the given diagram $K$,
 and the sum is taken over all possible colorings.
 The values of the state sum
 are  taken to be in  the group ring $\Z[A]$ where $A$ is the coefficient 
 group  written multiplicatively. 
 The state sum depends on the choice of $2$-cocycle $f$. 
 The cocycle invariant $\Phi_{\rm YB} (K)$ does not depend on 
 the choice of a diagram for $K$, %%%defining 
and thus is  %%
an invariant for virtual knots and links.
 Note that the image  of $ \Phi_{\rm YB} (K) $
 under the map $\Z[A] \rightarrow \Z$ sending all elements of $A$  to $1$
 is equal to the number of colorings  $| {\rm Col}_X(K)|$.

%%%
%% For the rest of the section, w
We investigate the 
 cocycle invariant for $G=\Z_n=A$ 
 with the %%%%
 $2$-cocycle $f(x,y)=x+y$. 

 \begin{theorem}\label{trivthm}
For the cocycle $f(x,y)=x+y$, suppose that $n$  is odd or $n=0.$ If a virtual knot, $K$, is checkboard coloarable, then %%%%$\Phi_{\rm YB} (K)
%%%%$ is a positive integer. That is $
%%=  %%%
$\sum_\tau B(\tau,{\cal C})=1$ for any coloring ${\cal C}.$ 
 \end{theorem}
 {\it Proof.\/} 
 Let $K$ be an 
oriented 
 virtual link diagram that has a  (mod $2$)-Alexander numbering,
 and a coloring ${\mathcal C}$ be given by a cyclic Wada biquandle $\Z_n$ with 
 $n$ odd. Consider the corresponding Fox coloring ${\mathcal C}_F$
 that is given by Lemma~\ref{Foxlemma}.
 Let $L: {\mathcal E} \rightarrow \Z_2$ be a (mod $2$)-Alexander numbering.
 Let ${\mathcal A}$ be the set of arcs in $K$.  (An arc is 
%%%neither broken 
broken neither at an over-crossing, nor at a virtual crossing). 
 To the initial and terminal points $s(\alpha)$, $t(\alpha)$, respectively, 
 of every arc $\alpha$, we make the following
 assignment $g$ of elements of $\Z_n$: 
 If  the Fox coloring ${\mathcal C}_F (\alpha)=x \in \Z_n$, then $g(s(\alpha))=- x$
 and $g(t(\alpha))= x$. 
 Since these have opposite signs, the sum over all assignments on 
 initial and terminal points of all arcs is zero.
 On the other hand, we show that the sum is equal to twice the weight at each crossing,
 and the result follows.
 Consider a positive crossing formed by an incoming arc $\alpha$ and the over-arc
 $\beta$ with ${\mathcal C} (\alpha)=x$ and ${\mathcal C} (\beta)=y$.
 Let $\gamma$ be the outgoing under-arc, so that ${\mathcal C} (\gamma)=x+2y$.
 Then $g(t(\alpha))= (-1)^{L(\alpha)} x$ and 
 $g(s(\gamma))=- (-1)^{L(\gamma)} (x+2y)$, 
 so the sum is $ (-1)^{L(\alpha)} 2(x+y)$ since $L(\gamma)=L(\alpha)+1$ 
which is the desired value (twice the weight of the crossing). 
The weight at a negative crossing is checked similarly.
 $\Box$

 As an application, let $VT(2, k)$ be an oriented 
  virtual $(2,k)$-torus knot 
or link 
represented by
 the closure of virtual $2$-braid $(\sigma_1)^k v_1$, where $\sigma_1$ 
 and $v_1$ denote the standard and virtual braid generators, respectively.
 The orientations are chosen to be downward.

 \begin{proposition}\label{torusprop}
 For any integer $k\neq 0$, $VT(2,k)$ is not checkerboard colorable.
 \end{proposition}
 {\it Proof.\/}
 Assume $k$ to be positive, as the negative case is similar.
 Assume that $VT(2,k)$ is almost classical (mod $2$) 
to  
derive a contradiction
 to Theorem~\ref{trivthm}. 
 Suppose the top left and right arcs receive colors 
$(x,y)$, $x, y\in \Z_n$. 
 Then below $(\sigma_1)^k$, the colors are 
 $(- ((k-1)x+ky), kx + (k+1)y)$. After the virtual crossing $v_1$,
 the bottom colors are $(kx + (k+1)y, - ((k-1)x+ky) )$, which must be equal to 
 $(x,y)$ for this to color the closure. 
 Thus $(x,y)$ colors if and only if $(k-1)x+(k+1)y \equiv 0$ \ (mod $n$).
  For a given $k$, choose $n$ such that ${\rm gcd}(k-1, n)=1$, 
  then there are $n$ solutions $(x,y)$ to the above equation.
  This implies $| {\rm Col}_X(K)|=n$, where $X=\Z_n$.
  One solution is $(x,y)=(k+1, 1-k)$, and the weight at the top crossing 
  (therefore at all the other crossings) is $2$.
  Hence the contribution of this particular coloring to the invariant is
  $2(k-1)$, which is not zero in $\Z_n$, giving a non-trivial contribution
  to the invariant. 
  $\Box$

 \begin{remark} {\rm
 Theorem~\ref{trivthm} can be stated more generally for 
 possibly infinite quandles:  
If $K$ is checkerboard  
 colorable, then for any coloring, the   product $\prod_{\tau}  B( \tau, {\cal C})$
 is the identity element of $A$. 
 Then  
in Proposition~\ref{torusprop}, the biquandle $\Z$ can be used,
 as the color vector $(k+1, 1-k) \in \Z \times \Z$ extends as well, 
 and gives the non-trivial product $2k$. Note that $|2k|$ is the span 
 of this coloring (Remark~\ref{spanrem}). 
 Thus we conjecture that the non-zero minimum of the product 
 $\prod_{\tau}  B( \tau, {\cal C})$ for the biquandle $\Z$ gives a lower bound
 of the span, and that the span of $VT(2,k)$ is $|2k|$.
 } \end{remark}

%MS% updated 
%%thanks

 \end{document}